\documentstyle[11pt]{article}

\newtheorem{definition}{Definition}[section]

\newtheorem{theorem}{Theorem}[section]

\begin{document}

\begin{center}

{\large\bf TWO-PARAMETER DEFORMATIONS OF LOOP ALGEBRAS 
AND SUPERALGEBRAS}

\vskip 12pt

{\Large Valeriy N. Tolstoy}
\vskip 3pt
{\large \it Institute of Nuclear Physics, Moscow State University\\
119899 Moscow\&Russia} {\it (e-mail: tolstoy@anna19.npi.msu.su)}

\end{center}

\begin{abstract}
We discuss two-parameter deformations of an universal enveloping 
algebra $U(g[u])$ of a polynomial loop algebra $g[u]$, 
where $g$ is a finite-dimensional complex simple  Lie algebra
(or superalgebra). These deformations are Hopf algebras. 
One deformation called Drinfeldian is a quantization of $U(g[u])$ 
in the direction of a classical r-matrix which is a sum of 
the simplest rational and trigonometric r-matrices. Another 
deformation (discussed only for the case $g=sl_{2}$) is 
a twisting of the usual Yangian $Y_{\eta}(sl_{2})$.
\end{abstract}
\section{Introduction}
As it is well known, an universal enveloping algebra
$U(g[u])$ of a polynomial loop (current) Lie algebra $g[u]$ , 
where $g$ is a finite-dimensional complex simple  Lie algebra,
admits two type deformations: a trigonometric deformation
$U_q(g[u])$ and a rational deformation or 
Yangian $Y_{\eta}(g)$ \cite{D}.
(In the case $g=sl_n$ there also exists an elliptic quantum
deformation of $U(sl_n[u])$, 
which is not discussed in this report).
The algebras $U_q(g[u])$, and $Y_{\eta}(g)$ are quantizations
of $U(g[u])$ in the direction of the simplest trigonometric
and rational solutions of the classical Yang-Baxter equation 
over $g$, respectively. These deformations are one-parameter.
It turns out that $U(g[u])$ also admits two-parameter
deformations. Here we discuss two type of such deformations 
which are Hopf algebras.

\noindent
A Hopf algebra of the first type called 
the rational-trigonometric quantum algebra or the Drinfeldian
$D_{q\eta}(g)$ \cite{T} is a quantization of $U(g[u])$ 
in the direction of a classical r-matrix which is a sum of 
the simplest rational and trigonometric r-matrices.
The  Drinfeldian $D_{q\eta}(g)$ contains $U_{q}(g)$ as a Hopf
subalgebra, and  $U_{q}(g[u])$ and $Y_{\eta}(g)$ are its limit
quantum algebras when the deformation parameters of $D_{q\eta}(g)$ 
$\eta$ goes to $0$ and $q$ goes to $1$, respectively. 
These results are easy generalized to a supercase, i.e. when 
$g$ is a finite-dimensional contragredient simple superalgebra.

\noindent  
A Hopf algebra of the second type discussed only for the case 
$g=sl_{2}$ is obtained by twisting of the usual Yangian 
$Y_{\eta}(sl_{2})$. The twisted Yangian $Y_{\eta\zeta}(sl_{2})$ 
is a quantization $U(sl_{2}[u])$ in the direction of a classical 
r-matrix 
$r(u,v)=\eta{\bf c}_{2}/(u-v)+\zeta h_{\alpha}\wedge e_{-\alpha}$,
where ${\bf c}_{2}$ is the $sl_{2}$ Casimir element. Detailed 
describtion of $Y_{\eta\zeta}(sl_{2})$ is given in \cite{KST}.
\section{Drinfeldian $D_{q\eta}(g)$}

Let $g$ be a finite-dimensional complex  simple Lie algebra
of a rank $r$ with a standard Cartan matrix $A=(a_{ij})_{i,j=1}^r$,
with a system of simple roots  $\Pi:= \{\alpha_1,\ldots, a_r\}$, 
and with a maximal positive root $\theta$.
Let $U_{q}(g)$ be a standard q-deformation of
the universal enveloping algebra $U(g)$ with 
Chevalley generators $k_{\alpha_i}^{\pm 1}$, $e_{\pm\alpha_i}$
$(i=1,2,\ldots, r)$ and with the defining relations 
\vspace{-3pt}
\begin{equation}
[k_{\alpha_i},k_{\alpha_j}]=0~,\qquad
k_{\alpha_i}e_{\pm\alpha_j}k^{-1}_{\alpha_i}=
q^{\pm(\alpha_i,\alpha_j)}e_{\pm\alpha_j}~,
\label{D1}
\end{equation}
\begin{equation}
[e_{\alpha_i},e_{-\alpha_i}]=
\frac{k_{\alpha_i}-k_{\alpha_i}^{-1}}{q-q^{-1}}~,\qquad
({\rm ad}_{q}e_{\pm\alpha_{i}})^{1-a_{ij}} e_{\pm\alpha_{j}}=0
\quad\;{\rm for}\,\, i\neq j~, 
\label{D2}
\end{equation}
\begin{definition}
The Drinfeldian $D_{q\eta}(g)$ is generated as an associative 
algebra over $C\!\!\!\!I\,[[\eta]]$ by the algebra $U_{q}(g)$ and 
the elements $\xi_{\delta-\theta}$, $k_{\delta}^{\pm 1}$ with 
the relations:
\vspace{-3pt}
\begin{equation}
[k_{\delta}^{\pm 1},{\rm everything}]=0~\qquad
k_{\alpha_i}\xi_{\delta-\theta}k^{-1}_{\alpha_i}=
q^{-(\alpha_i,\theta)}\xi_{\delta-\theta}~,
\label{D3}
\end{equation}
\begin{equation}
[e_{-\alpha_i},\xi_{\delta-\theta}]=
a\,[e_{-\alpha_i},\tilde{e}_{-\theta}],\qquad
({\rm ad}_{q}e_{\alpha_i})^{n_{i0}}\xi_{\delta-\theta}=
a\,({\rm ad}_{q}e_{\alpha_i})^{n_{i0}}\tilde{e}_{-\theta}
\label{D4}
\end{equation} 
for $n_{i0}=1+2(\alpha_i,\theta)/(\alpha_i,\alpha_i)$~, and
\begin{eqnarray}
[[e_{\alpha_i},\xi_{\delta-\theta}]_{q},\xi_{\delta-\theta}]_{q}
&\!\!\!\!\!=&\!\!\!\!\!
-a^{2}[[e_{\alpha_i},\tilde{e}_{-\theta}]_{q},\tilde{e}_{-\theta}]_{q}
\nonumber
\\
&&\!\!\!\!
+a\,[[e_{\alpha_i},\tilde{e}_{-\theta}]_{q},\xi_{\delta-\theta}]_{q}+
a\,[[e_{\alpha_i},\xi_{\delta-\theta}]_{q},\tilde{e}_{-\theta}]_{q}
\end{eqnarray}
\label{D5}
for $g\ne sl_2$ and $(\alpha_i,\theta)\ne 0$,
\begin{eqnarray}
&&\!\!\!\!\!\!\!\!\!\!\!\!\!
[[[e_{\alpha},\xi_{\delta-\alpha}]_{q},
\xi_{\delta-\alpha}]_{q},\xi_{\delta-\alpha}]_{q}=a^{3}
[[[e_{\alpha},\tilde{e}_{-\alpha}]_{q},\tilde{e}_{-\alpha}]_{q},
\tilde{e}_{-\alpha}]_{q}\qquad
\label{D6}
\nonumber
\\
&&-a^{2}[[[e_{\alpha},\tilde{e}_{-\alpha}]_{q},\tilde{e}_{-\alpha}]_{q},
\xi_{\delta-\alpha}]_{q}-a^{2}[[[e_{\alpha},\tilde{e}_{-\alpha}]_{q},
\xi_{\delta-\alpha}]_{q},\tilde{e}_{-\alpha}]_{q}
\nonumber\\
&&-a^{2}[[[e_{\alpha},\xi_{\delta-\alpha}]_{q},\tilde{e}_{-\alpha}]_{q},
\tilde{e}_{-\alpha}]_{q}+a\,[[[e_{\alpha},\tilde{e}_{-\alpha}]_{q},
\xi_{\delta-\alpha}]_{q},\xi_{\delta-\alpha}]_{q}
\nonumber
\\
&&+a\,[[[e_{\alpha},\xi_{\delta-\alpha}]_{q},
\tilde{e}_{-\alpha}]_{q},\xi_{\delta-\alpha}]_{q}+a\,[[[e_{\alpha},
\xi_{\delta-\alpha}]_{q},\xi_{\delta-\alpha}]_{q},
\tilde{e}_{-\alpha}]_{q}
\end{eqnarray}
for $g=sl_2$. The Hopf structure of $D_{q\eta}(g)$ is defined by 
the formulas $\Delta_{q\eta}(x)=\Delta_{q}(x)$, 
$S_{q\eta}(x)=S_{q}(x)$ ($x\in U_{q}(g)$) and it is the same for the
elements $k_{\delta}^{\pm}$ and  $k_{\alpha_i}$.
The comultiplication and the antipode of $\xi_{\delta-\alpha}$  are 
given by
\begin{eqnarray} 
\Delta_{q\eta}(\xi_{\delta-\theta})\!\!\!&=&\!\!\!
\xi_{\delta-\theta}\otimes 1+k_{\delta-\theta}^{-1}
\otimes \xi_{\delta-\theta}
\nonumber
\\
&&\!\!\!+ a\left(\Delta_{q}(\tilde{e}_{-\theta}) 
-\tilde{e}_{-\theta}\otimes 1-
k_{\delta-\theta}^{-1}\otimes\tilde{e}_{-\theta}\right),
\label{D7}
\end{eqnarray} 
\begin{equation}
S_{q\eta}(\xi_{\delta-\theta})=
-k_{\delta-\theta}\xi_{\delta-\theta}
+a\left(S_{q}(\tilde{e}_{-\theta})+
k_{\delta-\theta}\tilde{e}_{-\theta}\right)~,
\label{D8}
\end{equation}
where $a:=\eta/(q-q^{-1})$, 
$({\rm ad}_{q}e_{\beta})e_{\gamma}=[e_{\beta},e_{\gamma}]_{q}$, 
and the vector  $\tilde{e}_{-\theta}$ is any $U_{q}(g)$ 
element of the  weight $-\theta$, such that 
$g\ni\lim_{q\to1}\tilde{e}_{-\theta}\neq 0$.
\end{definition}
The right-hand sides of the relations (\ref{D4})-(\ref{D8}) are 
nonsingular at $q=1$.
\begin{theorem}
(i) The Drinfeldian $D_{q\eta}(g)$ is a two-parameter quantization
of $U(\overline{g[u]})$, where $\overline{g[u]}$ is a central extension 
of $g[u]$, in the direction of a classical r-matrix which
is a sum of the simplest rational and trigonometric r-matrices.
\\
\noindent
(ii) The Hopf algebra $D_{q=1,\eta}(g)$
is isomorphic to the Yangian $Y_{\eta}'(g)$ (with a central element). 
Moreover, $D_{q\eta=0}(g)=U_q(\overline{g[u]})$.
\label{DT1}
\end{theorem}

In the supercase, i.e. when $g$ is a simple finite-dimensional
contragredient Lie superalgebra all the commutators and 
the q-commutators are replaced by the supercommutators and 
the q-supercommutators. Moreover we have to add some additional
Serre relations if they exist. 

\section{Twisted Yangian $Y_{\eta\zeta}(sl_2)$}

In the case $sl_{2}$ from (\ref{D3})-(\ref{D8}) at $q=1$ we can
obtain that the Yangian $Y_{\eta}(sl_2)$ is generated by 
the $sl_2$ elements $h_{\alpha},\;e_{\pm\alpha}$ and the element 
$\xi_{\delta-\alpha}$ with the relations:
\begin{equation}
[h_{\alpha},\xi_{\delta-\alpha}]=-2\xi_{\delta-\alpha}~,\quad
[e_{-\alpha},\xi_{\delta-\alpha}]=\eta e_{-\alpha}^2~,
\label{TY2}
\end{equation}
\begin{equation}
[e_{\alpha},[e_{\alpha},[e_{\alpha},\xi_{\delta-\alpha}]]]=
6\eta e_{\alpha}^{2}~,\quad
[[[e_{\alpha},\xi_{\delta-\alpha}],\xi_{\delta-\alpha}],\xi_{\delta-\alpha}]=
6\eta \xi_{\delta-\alpha}^{2}~.
\label{TY3}
\end{equation}
\begin{equation}
\Delta(\xi_{\delta-\alpha})=\xi_{\delta-\alpha}\otimes 1 + 
1\otimes \xi_{\delta-\alpha}+
\eta e_{-\alpha}\otimes h_{\alpha}~,\;\; 
S(\xi_{\delta-\alpha})=-\xi_{\delta-\alpha}+\eta e_{-\alpha}h_{\alpha}~,
\label{TY5}
\end{equation}
where we put $(\alpha,\alpha)=2$.

Using the twisting element $F=\sum_{k\geq 0}(\zeta^{k}/k!)
\big(\prod_{i=0}^{k-1}(h_{\alpha}+2i)\big)\otimes e_{-\alpha}^k~$ 
one can calculate the new coproduct $\Delta^{(F)}(x):=F\Delta(x)F^{-1}$ 
and antipode $S^{(F)}:=uS(x)u^{-1}$ ($x\in Y_{\eta}(sl_2)$) where 
$u:=\sum_{k\geq 0}((-\zeta)^k/k!)
\big(\prod_{i=0}^{k-1}(h_{\alpha}+2i)\big)e_{-\alpha}^k$~. 
The result is the twisted Yangian $Y_{\eta\zeta}(sl_2)$.
It is not difficult to show that $Y_{\eta\zeta}(sl_2)$ is
a quantization of $U(sl_{2}[u])$ with the classical r-matrix 
$r(u,v)=\eta{\bf c}_{2}/(u-v)+\zeta h_{\alpha}\wedge e_{-\alpha}$,
where ${\bf c}_{2}$ is the $sl_{2}$ Casimir element.
See \cite{KST} for details.

\section*{Acknowledgments}
The author is grateful to the Organizing Committee 
for a partial support of his participation at WigSym5.
The work was supported by  the Russian  Foundation 
for Fundamental Research, grant No. 96-01-01421.


\end{document}